\documentclass{article}
\usepackage{amsmath,amsthm,amsopn,amstext,amscd,amsfonts,amssymb}
\usepackage{natbib}
\usepackage{verbatim}

\def\diag{\mathop{\rm diag}\nolimits}

\def\eig{\mathop{\rm ch}\nolimits}

\def\build#1#2#3{\mathrel{\mathop{#1}\limits^{#2}_{#3}}}

\def \Prod#1#2{\mbox{$\displaystyle \build{\prod}{#2}{#1}$}}

\newcommand {\findemo}{\hfill \, $\Box$ \,\\[2ex]}

\renewenvironment{abstract}
                 {\vspace{6pt}
                  \begin{center}
                  \begin{minipage}{5in}
                  \centerline{\textbf{Abstract}}
                  \noindent\ignorespaces
                 }
                 {\end{minipage}\end{center}}
\newtheorem{thm}{\textbf{Theorem}}[section]

\newtheorem{lem}{\textbf{Lemma}}[section]
\theoremstyle{definition}

\newtheorem{rem}{\textbf{Remark}}[section]

\title{\huge \textbf{Doubly singular matrix variate beta type I and II and
        singular inverted matricvariate $t$ distributions}}
\author{
  \textbf{Jos\'e A. D\'{\i}az-Garc\'{\i}a} \thanks{Corresponding author\newline
    }\\
  Department of Statistics and Computation \\
  25350 Buenavista, Saltillo, Coahuila, Mexico \\
  E-mail: jadiaz@uaaan.mx \\[2ex]
  \textbf{Ram\'on Guti\'errez J\'aimez} \\
  Department of Statistics and O.R. \\
  University of Granada \\
  Granada 18071, Spain \\
  E-mail: rgjaimez@ugr.es\\
}
\date{}
\begin{document}
\maketitle

\begin{abstract}
In this paper, the densities of the doubly singular beta type I and II distributions
are found, and the joint densities of their corresponding nonzero eigenvalues are
provided. As a consequence, the density function of a singular inverted matricvariate
$t$ distribution is obtained.
\end{abstract}

\section{Introduction}\label{Intro}

Let $\mathbf{A}$ and $\mathbf{B}$ be independent Wishart matrices or, alternatively,
let  $\mathbf{B} = \mathbf{Y}\mathbf{Y}'$ where $\mathbf{Y}$ has a matrix variate
normal distribution. Then the matrix variate beta type I distribution is defined in
the literature as
\begin{equation}\label{defbetaI}
  \mathbf{U} = \left\{
  \begin{array}{l}
    (\mathbf{A} + \mathbf{B})^{-1/2}\mathbf{B} (\mathbf{A} + \mathbf{B})^{-1/2}, \\
    \mathbf{B}^{1/2}(\mathbf{A} + \mathbf{B})^{-1} \mathbf{B}^{1/2}, \\
    \mathbf{Y}'(\mathbf{A} + \mathbf{B})^{-1}\mathbf{Y}.
  \end{array}
  \right.
\end{equation}
Analogously, the matrix variate beta type II distribution is defined as
\begin{equation}\label{defbetaII}
  \mathbf{F} = \left\{
  \begin{array}{l}
    \mathbf{A}^{-1/2}\mathbf{BA}^{-1/2}, \\
    \mathbf{B}^{1/2}\mathbf{A}^{-1} \mathbf{B}^{1/2}, \\
    \mathbf{Y}' \mathbf{A}^{-1}\mathbf{Y}.
  \end{array}
  \right.
\end{equation}
These can be classified as central, noncentral or doubly noncentral, depending on
whether $\mathbf{A}$ and $\mathbf{B}$ are central, $\mathbf{B}$ is noncentral or
$\mathbf{A}$ and $\mathbf{B}$ are noncentral, respectively; see \citet{dggj:07,
dggj:08b}. In addition, such a distribution can be classified as nonsingular,
singular or doubly singular, when $\mathbf{A}$ and $\mathbf{B}$ are nonsingular,
$\mathbf{B}$ is singular or $\mathbf{A}$ and $\mathbf{B}$ are singular, respectively,
see \citet{dggj:08a}. Under definitions (\ref{defbetaI}) and (\ref{defbetaII}) and
their classifications, mentioned above, the matrix variate beta type I and II
distributions have been studied by different authors; see \citet{or:64},
\citet{di:67}, \citet{m:70}, \citet{k:70}, \citet{sk:79}, \citet{mh:82},
\citet{u:94}, \citet{c:96}, \citet{dg:97, dggj:06, dggj:07, dggj:08a, dggj:08b},
among many others. However, the density functions of doubly singular beta type I and
II distributions have not been studied. In such cases the inverses that appear in
definitions (\ref{defbetaII}) and (\ref{defbetaI}) must be replaced by the
Moore-Penrose inverse. The doubly singular beta type I distribution was briefly
considered by \citet{m:70}, but this distribution has been the object of considerable
recent interest. Observe that this distribution appears in a natural way when the
number of observations is smaller than the dimension in multivariate analysis of
variance; see \citet{s:07}.

Furthermore, observe that if $\mathbf{Z}$ is any random matrix  $m \times r$, it is
generically called as matrix variate if the kernel of its density function is in
terms of the trace operator. But if the kernel of the density function is in terms of
the determinant alone or in terms of both, then the determinant and the trace
operator, $\mathbf{Z}$, are generically called matricvariate, see \citet{di:67} and
\citet{dggj:08c}. Alternatively, when the random matrix $\mathbf{A}: m \times m$ is
symmetric by definition (e.g. Wishart and beta matrices, etc.), they are generically
referred to as matrix variate.

In Section \ref{sec2} of the present study, we review singular matricvariate $t$
distributions and determine the distribution of a linear transformation and the joint
density function of its nonzero singular values. In Section \ref{sec3}, an expression
is provided for the density function of the doubly singular matrix variate beta type
II distribution, and the joint density function of its nonzero eigenvalues is
determined.

Similar results are provided for the doubly singular matrix variate beta type I
distribution; see Section \ref{sec4}. Finally, in Section \ref{sec5}, we study the
inverted matricvariate $t$ distribution, also termed the matricvariate Pearson type
II distribution.

\section{Preliminary results}\label{sec2}

Let ${\mathcal L}_{r,m}(q)$ be the linear space of all $m \times r$ real matrices of
rank $q \leq \min(m,r)$ and let ${\mathcal L}_{r,m}^{+}(q)$ be the linear space of
all $m \times r$ real matrices of rank $q \leq \min(m,r)$, with $q$ distinct singular
values. The set of matrices $\mathbf{H}_{1} \in {\mathcal L}_{r,m}$ such that
$\mathbf{H}'_{1}\mathbf{H}_{1} = I_{r}$ is a manifold denoted as ${\mathcal
V}_{r,m}$, termed the Stiefel manifold. In particular, ${\mathcal V}_{m,m}$ is the
group of orthogonal matrices ${\mathcal O}(m)$. The invariant measure on a Stiefel
manifold is given by the differential form
$$
    (\mathbf{H}'_{1}d\mathbf{H}_{1})\equiv \bigwedge_{i=1}^{r} \bigwedge_{j =i+1}^{m}
        \mathbf{h}'_{j}d\mathbf{h}_{i}
$$
written in terms of the exterior product ($\bigwedge$), where we choose an $m \times
(m-r)$ matrix $\mathbf{H}_{2}$ such that $\mathbf{H}$ is an $m \times m$ orthogonal
matrix, with $\mathbf{H} = (\mathbf{H}_{1} \vdots \mathbf{H}_{2})$ and where
$d\mathbf{h}$ is an $m \times 1$ vector of differentials; see \citet[Section
2.1.4]{mh:82}. Moreover
\begin{equation}\label{vol}
    \int_{\mathbf{H}_{1}\in {\mathcal V}_{r,m}}(\mathbf{H}'_{1}d\mathbf{H}_{1}) = \frac{2^{r}
    \pi^{mr/2}}{\Gamma_{r}[m/2]} \quad \mbox{and} \quad \int_{\mathbf{H} \in {\mathcal O}(m)}
    (\mathbf{H}'d\mathbf{H}) = \frac{2^{m}  \pi^{m^{2}/2}}{\Gamma_{m}[m/2]}
\end{equation}

Denote by ${\mathcal S}_{m}$, the homogeneous space of $m \times m$ positive definite
symmetric matrices; and by ${\mathcal S}_{m}(q)$, the $(mq - q(q - 1)/2)$-dimensional
manifold of rank $q$ positive semidefinite $m \times m$ symmetric matrices; and by
${\mathcal S}_{m}^{+}(q)$, the $(mq - q(q - 1)/2)$-dimensional manifold of rank $q$
positive semidefinite $m \times m$ symmetric matrices with $q$ distinct positive
eigenvalues. Assume $\mathbf{A} \in {\mathcal S}_{m}$; then $\eig_{i}(\mathbf{A})$
denotes the $i$-th eigenvalue of the matrix $\mathbf{A}$. Moreover, let
$\mathbf{A}^{+}$ and $\mathbf{A}^{-}$ be the Moore-Penrose inverse and any symmetric
generalized inverse of $\mathbf{A}$, respectively; see \citet{r:73}. Finally,
$\mathbf{A}^{1/2}$ is termed a non-negative definite square root of $\mathbf{A}$ if
$\mathbf{A}^{1/2}$ is such that $\mathbf{A} = \mathbf{A}^{1/2}\mathbf{A}^{1/2}$.

Using the notation of \citet{dgm:97} for the singular matrix variate normal and
Wishart and pseudo-Wishart distributions, from \citet{dggj:08c} we have:

\begin{lem}\label{lemma1}
Let $\mathbf{T} \in {\mathcal L}_{r,m}^{+}(r_{_{\mathbf{\Xi}}})$ be the random matrix
$$
  \mathbf{T} = \left (\mathbf{A}^{1/2} \right)^{+}\mathbf{Y} + \mathbf{\mu}
$$
where $\mathbf{A}^{1/2}\mathbf{A}^{1/2} = \mathbf{A} \in {\mathcal S}_{m}^{+}(q) \sim
\mathcal{W}_{m}^{q}\left(n,\mathbf{\Theta} \right)$, $\mathbf{\Theta} \in {\mathcal
S}_{m}(r_{_{\mathbf{\Theta}}})$ with $r_{_{\mathbf{\Theta}}} \leq m$ and $q =
\min(m,n)$, independent of  $\mathbf{Y}\in {\mathcal
L}_{r,m}^{+}(r_{_{\mathbf{\Xi}}}) \sim \mathcal{N}_{m \times r}^{m,
r_{_{\mathbf{\Xi}}}}(\mathbf{0}, \mathbf{I}_{m} \otimes \mathbf{\Xi})$, $\mathbf{\Xi}
\in {\mathcal S}_{r}(r_{_{\mathbf{\Xi}}})$ with $ r_{_{\mathbf{\Xi}}} \leq r$ and $m
\geq q \geq r_{_{\mathbf{\Xi}}} > 0$. Then the singular matricvariate $\mathbf{T}$
has the density
\begin{equation}\label{mcT}
   \hspace{-1cm}
    \frac{c(m,n,q,q_{1},r_{_{\mathbf{\Theta}}}, r_{_{\mathbf{\Xi}}}, r_{_{\alpha}})}
    {\Prod{i=1}{r_{_{\mathbf{\Xi}}}}\eig_{i}(\mathbf{\Xi})^{m/2} \ \Prod{j=1}{r_{_{\mathbf{\Theta}}}}
    \eig_{j}(\mathbf{\Theta})^{n/2}}\prod_{l = 1}^{r_{_{\alpha}}}\eig_{l}\left[\mathbf{\Theta}^{-} +
    (\mathbf{T} - \mathbf{\mu}) \mathbf{\Xi}^{-}(\mathbf{T} - \mathbf{\mu})'
    \right]^{-\left(n+r_{_{\mathbf{\Xi}}}\right)/2}(d\mathbf{T}),
\end{equation}
with
\begin{equation}\label{cte}
\hspace{-.5cm}
    c(m,n,q,q_{1},r_{_{\mathbf{\Theta}}}, r_{_{\mathbf{\Xi}}}, r_{_{\alpha}})=
    \frac{\pi^{n(q-r_{_{\mathbf{\Theta}}})/2-(n+r_{_{\mathbf{\Xi}}})
    (q_{1}-r_{_{\alpha}})/2-m r_{_{\mathbf{\Xi}}}/2}
    \Gamma_{q_{1}}[(n+r_{_{\mathbf{\Xi}}})/2]}{2^{(m r_{_{\mathbf{\Xi}}}+
    n r_{_{\mathbf{\Theta}}})/2-(n+r_{_{\mathbf{\Xi}}})r_{_{\alpha}}/2} \Gamma_{q}[n/2]},
\end{equation}
where $q_{1} = \min(m,n+r_{_{\mathbf{\Xi}}})$, $r_{_{\alpha}} =
\hbox{rank}\left[\mathbf{\Theta}^{-} + (\mathbf{T} - \mathbf{\mu})
\mathbf{\Xi}^{-}(\mathbf{T} - \mathbf{\mu})'\right] \leq m$, and $(d\mathbf{T})$
denotes the Hausdorff measure.
\end{lem}

In particular observe that if $\mathbf{\Xi} = \mathbf{I}_{r}$, i.e.
$r_{_{\mathbf{\Xi}}} = r$, $r_{\mathbf{\Theta}} = r_{\alpha} = m$, $\mathbf{\Theta} =
\mathbf{I}_{m}$, $q = n$ and $q_{1} = n +r$, we have
\begin{equation}\label{Tdist}
    dF_{\mathbf{T}}(\mathbf{T}) = \pi^{-r(r+2n)/2}
    \frac{\Gamma_{n+r}[(n+r)/2]}{\Gamma_{n}[n/2]} |\mathbf{I}_{m}+
    (\mathbf{T}-\mathbf{\mu})(\mathbf{T}- \mathbf{\mu})'|^{-(n+r)/2}(d\mathbf{T})
\end{equation}
where $\mathbf{T} \in {\mathcal L}_{r,m}^{+}(r)$ and now $(d\mathbf{T})$ denotes the
Lebesgue measure.

\begin{thm}\label{teo2}
Under the condition of Lemma \ref{lemma1}, let $\mathbf{\mu} = \mathbf{0}$ and
$\mathbf{X} = \mathbf{TC}^{+'} \in {\mathcal
L}_{r_{_{\mathbf{\Xi}}},m}^{+}(r_{_{\mathbf{\Xi}}})$, where $\mathbf{\Xi} =
\mathbf{CC}'$ such that $\mathbf{C} \in {\mathcal
L}_{r_{_{\mathbf{\Xi}}},r}^{+}(r_{_{\mathbf{\Xi}}})$. Then the density function of
$\mathbf{X}$ is
\begin{equation}\label{mcTX}
    \frac{c(m,n,q,q_{1},r_{_{\mathbf{\Theta}}}, r_{_{\mathbf{\Xi}}}, r_{_{\alpha}})}
    {\Prod{j=1}{r_{_{\mathbf{\Theta}}}} \eig_{j}(\mathbf{\Theta})^{n/2}}\prod_{l = 1}^{r_{_{\alpha}}}
    \eig_{l}\left(\mathbf{\Theta}^{-} +
    \mathbf{X}\mathbf{X}'\right)^{-\left(n+r_{_{\mathbf{\Xi}}}\right)/2}(d\mathbf{X}),
\end{equation}
where $q_{1} = \min(m,n+r_{_{\mathbf{\Xi}}})$, $r_{_{\alpha}} =
\hbox{rank}\left(\mathbf{\Theta}^{-} + \mathbf{X}\mathbf{X}'\right) \leq m$, and now
$(d\mathbf{X})$ denotes the Lebesgue measure.
\end{thm}
\textit{Proof.} Make the change of variables $\mathbf{T} = \mathbf{XC}^{'}$ as
described by \citet{dg:07}
\begin{equation}\label{jaco}
    (d\mathbf{T}) = \Prod{i=1}{r_{_{\mathbf{\Xi}}}}\eig_{i}(\mathbf{CC}')^{m/2} (d\mathbf{X})
                = \Prod{i=1}{r_{_{\mathbf{\Xi}}}}\eig_{i}(\mathbf{\Xi})^{m/2}
                  (d\mathbf{X}).
\end{equation}
Also, observe that
\begin{equation}\label{eqC}
    \mathbf{C}^{+} \mathbf{\Xi} \mathbf{C}^{+'} = \mathbf{C}^{+} (\mathbf{CC}') \mathbf{C}^{+'} =
  \mathbf{C}^{+} \mathbf{C}(\mathbf{C}^{+} \mathbf{C})' = \mathbf{C}^{+} \mathbf{C} =
  \mathbf{I}_{r_{_{\mathbf{\Xi}}}},
\end{equation}
since $\mathbf{C}$ and $\mathbf{C}^{+}\mathbf{C}$ have the same rank
$r_{_{\mathbf{\Xi}}}$ and $\mathbf{C}^{+}\mathbf{C}$ is of the order
$r_{_{\mathbf{\Xi}}} \times r_{_{\mathbf{\Xi}}}$. Finally, note that by (\ref{eqC})
\begin{equation}\label{eqC2}
    \mathbf{T} \mathbf{\Xi}^{-} \mathbf{T}' =\mathbf{XC'} \mathbf{\Xi}^{-} \mathbf{CX}'
  = \mathbf{X} (\mathbf{C}^{-}\mathbf{\Xi C}^{'-})^{-}\mathbf{X}' = \mathbf{X} (\mathbf{C}^{+}
  \mathbf{\Xi C}^{+'})^{-}\mathbf{X}' = \mathbf{X} \mathbf{X}'.
\end{equation}
By substituting (\ref{jaco}) and (\ref{eqC2}) in (\ref{mcT}) we obtain the desired
result. \findemo

Now let $\kappa_{i}$, $i = 1, \dots, r$ be the nonzero singular value of
$\mathbf{T}$. Then taking $\mathbf{\mu} = \mathbf{0}$ in (\ref{Tdist}) we have
\begin{thm}\label{teo1}
The joint density function of $\kappa_{1}, \dots,\kappa_{r}$, the singular values of
$\mathbf{T}$ is
\begin{equation}\label{svT}\hspace{-.5cm}
    \pi^{-r(2n-m)/2}\frac{2^{r} \ \Gamma_{n+r}[(n+r)/2]}{\Gamma_{n}[n/2] \Gamma_{r}[r/2]
    \Gamma_{r}[m/2]} \prod_{i = 1}^{r}\frac{\kappa_{i}^{m-r}}{(1+\kappa_{i}^{2})^{(n+r)/2}}
    \prod_{i<j}\left(\kappa_{i}^{2} -\kappa_{j}^{2}\right)\left(\bigwedge_{i=1}^{r} d\kappa_{i}\right)
\end{equation}
where $\kappa_{1}> \cdots >\kappa_{r}>0$.
\end{thm}
\textit{Proof.} Let $\mathbf{T} = \mathbf{Q}_{1} \mathbf{D}_{\mathbf{T}} \mathbf{P}'$
be the nonsingular part of the singular value decomposition (SVD), where
$\mathbf{Q}_{1} \in \mathcal{V}_{r,m}$, $\mathbf{P} \in \mathcal{O}(r)$ and
$\mathbf{D}_{\mathbf{T}} = \diag(\kappa_{1}, \dots,\kappa_{r})$, $\kappa_{1}> \cdots
>\kappa_{r}>0$. Then, from \citet{dgm:97}
$$
  (d\mathbf{T}) = 2^{-r}|\mathbf{D}_{\mathbf{T}}|^{m-r}\prod_{i<j}\left(\kappa_{i}^{2} -\kappa_{j}^{2}\right)
  (d\mathbf{D}_{\mathbf{T}})(\mathbf{Q}'_{1}d\mathbf{Q}_{1})(\mathbf{P}'d\mathbf{P}).
$$
The result is then obtained immediately from (\ref{Tdist}), noting that
$$
  |\mathbf{I}_{m} + \mathbf{Q}_{1} \mathbf{D}_{\mathbf{T}} \mathbf{P}'\mathbf{P}
    \mathbf{D}_{\mathbf{T}} \mathbf{Q}'_{1}| = |\mathbf{I}_{r} + \mathbf{D}_{\mathbf{T}}^{2}|
$$
and integrating over $ \mathbf{Q}_{1} \in \mathcal{V}_{r,m}$  and  $ \mathbf{P} \in
\mathcal{O}(r)$ using (\ref{vol}). \findemo

\section{Doubly singular beta type II distribution}\label{sec3}

\begin{thm}\label{teoBII1}
Let $\mathbf{A}^{1/2}\mathbf{A}^{1/2} = \mathbf{A} \in {\mathcal
S}_{m}^{+}(q_{_{\mathbf{A}}}) \sim
\mathcal{W}_{m}^{q_{_{\mathbf{A}}}}\left(n,\mathbf{\Theta} \right)$, $\mathbf{\Theta}
\in {\mathcal S}_{m}(r_{_{\mathbf{\Theta}}})$ with $r_{_{\mathbf{\Theta}}} \leq m$
and $q_{_{\mathbf{A}}} = \min(m,n)$, independent of  $\mathbf{Y}\in {\mathcal
L}_{r,m}^{+}(r_{_{\mathbf{\Xi}}}) \sim \mathcal{N}_{m \times r}^{m,
r_{_{\mathbf{\Xi}}}}(\mathbf{0}, \mathbf{I}_{m} \otimes \mathbf{\Xi})$, $\mathbf{\Xi}
\in {\mathcal S}_{r}(r_{_{\mathbf{\Xi}}})$ with $ r_{_{\mathbf{\Xi}}} \leq r$ and $m
\geq q \geq r_{_{\mathbf{\Xi}}} > 0$. The doubly singular matrix variate beta type II
distribution is defined as
$$
  \mathbf{F} = (\mathbf{A}^{1/2})^{+} (\mathbf{Y}\mathbf{\Xi}^{-}\mathbf{Y})
(\mathbf{A}^{1/2})^{+} = (\mathbf{A}^{1/2})^{+} \mathbf{B} (\mathbf{A}^{1/2})^{+}
$$
where $\mathbf{B} = (\mathbf{Y}\mathbf{\Xi}^{-}\mathbf{Y}) \in {\mathcal
S}_{m}^{+}(r_{_{\mathbf{\Xi}}}) \sim \mathcal{PW}_{m}^{r_{_{\mathbf{\Xi}}}}
\left(r_{_{\mathbf{\Xi}}}, \mathbf{I}_{m} \right)$. The density function of
$\mathbf{F}$ is
\begin{equation}\label{dfbetaII1} \hspace{-1cm}
    \frac{c(m,n,q,q_{1},r_{_{\mathbf{\Theta}}}, r_{_{\mathbf{\Xi}}}, r_{_{\alpha}})}
    {\pi^{-r_{_{\mathbf{\Xi}}}^{2}/2} \Gamma_{r_{_{\mathbf{\Xi}}}}[r_{_{\mathbf{\Xi}}}/2]
    \Prod{j=1}{r_{_{\mathbf{\Theta}}}} \eig_{j}
    (\mathbf{\Theta})^{n/2}} |\mathbf{D}_{\mathbf{F}}|^{(r_{_{\mathbf{\Xi}}}-m-1)/2}
    \prod_{l = 1}^{r_{_{\alpha}}}
    \eig_{l}\left(\mathbf{\Theta}^{-} +  \mathbf{F} \right)^{-\left(n+r_{_{\mathbf{\Xi}}}\right)/2}
    (d\mathbf{F}),
\end{equation}
where $\mathbf{F} = \mathbf{Q}_{1} \mathbf{D}_{\mathbf{F}} \mathbf{Q}'_{1}$ is the
nonsingular part of the spectral decomposition, with $\mathbf{Q}_{1} \in
\mathcal{V}_{r_{_{\mathbf{\Xi}}},m}$ and $\mathbf{D}_{\mathbf{F}} = \diag(\delta_{1},
\dots, \delta_{r_{_{\mathbf{\Xi}}}})$, $\delta_{1}> \cdots >
\delta_{r_{_{\mathbf{\Xi}}}}
>0$ and where $(d\mathbf{F})$ denotes the Hausdorff measure.
\end{thm}
\textit{Proof}. By (\ref{eqC2}), observe that
$$
  \mathbf{F} = (\mathbf{A}^{1/2})^{+} (\mathbf{Y}\mathbf{\Xi}^{-}\mathbf{Y})
    (\mathbf{A}^{1/2})^{+} = \mathbf{T} \mathbf{\Xi}^{-} \mathbf{T}' = \mathbf{X}\mathbf{X}'
$$
Let $\mathbf{X} = \mathbf{Q}_{1} \mathbf{D}_{\mathbf{X}} \mathbf{P}'$ be the
nonsingular part of the singular value decomposition (SVD), where $\mathbf{Q}_{1} \in
\mathcal{V}_{r_{_{\mathbf{\Xi}}},m}$, $\mathbf{P} \in \mathcal{O}(r)$ and
$\mathbf{D}_{\mathbf{X}} = \diag(\kappa_{1}, \dots,\kappa_{r_{_{\mathbf{\Xi}}}})$,
$\kappa_{1}> \cdots > \kappa_{r_{_{\mathbf{\Xi}}}}>0$. Then $\mathbf{F} =
\mathbf{X}\mathbf{X}' = \mathbf{Q}_{1} \mathbf{D}_{\mathbf{X}} \mathbf{P}'\mathbf{P}
\mathbf{D}_{\mathbf{X}} \mathbf{Q}_{1}' = \mathbf{Q}_{1} \mathbf{D}_{\mathbf{X}}^{2}
\mathbf{Q}_{1}' = \mathbf{Q}_{1} \mathbf{D}_{\mathbf{F}} \mathbf{Q}_{1}'$, with $
\mathbf{D}_{\mathbf{F}} =\diag(\delta_{1}, \dots, \delta_{r_{_{\mathbf{\Xi}}}})=
\mathbf{D}_{\mathbf{X}}^{2}$. Then, from \citet{dgm:97},
\begin{equation}\label{jacF}
    (d\mathbf{X}) = 2^{-r_{_{\mathbf{\Xi}}}}|\mathbf{D}_{\mathbf{F}}|^{(r_{_{\mathbf{\Xi}}} -m -1)/2}
        (d\mathbf{F})(\mathbf{P}'d\mathbf{P}).
\end{equation}
Substituting (\ref{jacF}) in (\ref{mcTX}) and integrating over $ \mathbf{P} \in
\mathcal{O}(r)$ using (\ref{vol}) gives the stated density function for $\mathbf{F}$.
\findemo

In particular, if $\mathbf{\Theta} = I_{m}$, $q_{1}= n + r_{_{\mathbf{\Xi}}} $ and $q
= n$, then $r_{\alpha} = m$ and
\begin{equation}\label{dfbetaII11}
    \pi^{-nr_{_{\mathbf{\Xi}}}}
    \frac{ \Gamma_{n+r_{_{\mathbf{\Xi}}}}[(n + r_{_{\mathbf{\Xi}}})/2]}
    {\Gamma_{n}[n/2]\Gamma_{r_{_{\mathbf{\Xi}}}}[r_{_{\mathbf{\Xi}}}/2]}
    |\mathbf{D}_{\mathbf{F}}|^{(r_{_{\mathbf{\Xi}}}-m-1)/2} |\mathbf{I}_{m} +
    \mathbf{F}|^{-\left(n+r_{_{\mathbf{\Xi}}}\right)/2}
    (d\mathbf{F}).
\end{equation}

An alternative way to obtain (\ref{dfbetaII11}) when $\mathbf{\Xi} = I_{r}$ and
therefore $r_{_{\mathbf{\Xi}}} = r$ is by using the third definition of $\mathbf{F}$
given in (\ref{defbetaII}). This denotes the matrix variate as
$\widetilde{\mathbf{F}}$.

\begin{thm}\label{teoBII2}
Let $\mathbf{A}^{1/2}\mathbf{A}^{1/2} = \mathbf{A} \in {\mathcal S}_{m}^{+}(n) \sim
\mathcal{PW}_{m}^{n}\left(n,\mathbf{I}_{m} \right)$, independent of  $\mathbf{Y}\in
{\mathcal L}_{r,m}^{+}(r) \sim \mathcal{N}_{m \times r}(\mathbf{0}, \mathbf{I}_{m}
\otimes \mathbf{I}_{r})$, and $m \geq n \geq r > 0$. The doubly singular matrix
variate beta type II distribution is defined as
$$
  \widetilde{\mathbf{F}} = \mathbf{Y}'\mathbf{A}^{+} \mathbf{Y},
$$
which has the density function
\begin{equation}\label{dfbetaII2} \pi^{-r(r+2n-m)/2}
    \frac{ \Gamma_{n+r}[(n + r)/2]}
    {\Gamma_{n}[n/2]\Gamma_{r}[m/2]}
    |\widetilde{\mathbf{F}}|^{(m- r-1)/2} |\mathbf{I}_{r} +
    \widetilde{\mathbf{F}}|^{-\left(n+r\right)/2}
    (d\widetilde{\mathbf{F}}),
\end{equation}
where $(d\widetilde{\mathbf{F}})$ denotes the Lebesgue measure.
\end{thm}
\textit{Proof}. Observe that
$$
  \widetilde{\mathbf{F}} = \mathbf{Y}'\mathbf{A}^{+} \mathbf{Y} = \mathbf{T}'\mathbf{T},
$$
where the density function of $\mathbf{T}$ is given by (\ref{Tdist}). Now as in the
proof of Theorem \ref{teo1}, let $\mathbf{T} = \mathbf{Q}_{1} \mathbf{D}_{\mathbf{T}}
\mathbf{P}'$; hence $\widetilde{\mathbf{F}} = \mathbf{T}'\mathbf{T} = \mathbf{P}
\mathbf{D}_{\mathbf{T}} \mathbf{Q}'_{1}\mathbf{Q}_{1} \mathbf{D}_{\mathbf{T}}
\mathbf{P}' = \mathbf{P}' \mathbf{D}_{\widetilde{\mathbf{F}}} \mathbf{P}'$, with $
\mathbf{D}_{\widetilde{\mathbf{F}}} = \mathbf{D}_{\mathbf{T}}^{2}$. Then by
\citet{dgm:97}
\begin{equation}\label{jacFt}
    (d\mathbf{T}) = 2^{-r}|\widetilde{\mathbf{F}}|^{m - r -1)/2}
        (d\widetilde{\mathbf{F}})(\mathbf{Q}'_{1}d\mathbf{Q}_{1}).
\end{equation}
The proof is completed by substituting (\ref{jacFt}) in (\ref{Tdist}) and integrating
over $ \mathbf{Q}_{1} \in \mathcal{V}_{r,m}$ using (\ref{vol}). \findemo

Note that Theorems \ref{teoBII1} and \ref{teoBII2} are generalizations of Theorems
10.4.1 and 10.4.4 in \citet{mh:82} for the central doubly singular case.

\begin{thm}\label{teochibetaII}
Let $\delta_{1}, \dots, \delta_{r}$  be the nonzero eigenvalues of $\mathbf{F}$ (with
density function (\ref{dfbetaII11}) and let $r_{_{\mathbf{\Xi}}} = r$) or equivalent
be the nonzero eigenvalues of $\widetilde{\mathbf{F}}$ (with density function
(\ref{dfbetaII2})). Then the joint density function of $\delta_{1}, \dots,
\delta_{r}$ is
\begin{equation}\label{chibetaII}\hspace{-.5cm}
    \pi^{-r(2n-m)/2}\frac{\Gamma_{n+r}[(n+r)/2]}{\Gamma_{n}[n/2] \Gamma_{r}[r/2]
    \Gamma_{r}[m/2]} \prod_{i = 1}^{r}\frac{\delta_{i}^{(m-r-1)/2}}{(1+\delta_{i})^{(n+r)/2}}
    \prod_{i<j}\left(\delta_{i} - \delta_{j}\right)\left(\bigwedge_{i=1}^{r} d\delta_{i}\right)
\end{equation}
where $\delta_{1} > \cdots > \delta_{r} > 0$.
\end{thm}
\textit{Proof}. This result can be obtained by any of the following methods:
\begin{description}
  \item[i)] By making the transformation $\mathbf{F} = \mathbf{Q}_{1}
    \mathbf{D}_{\mathbf{F}} \mathbf{Q}_{1}'$ in (\ref{dfbetaII11}) with $r_{_{\mathbf{\Xi}}} =
    r$ and integrating over $ \mathbf{Q}_{1} \in \mathcal{V}_{r,m}$ using
    (\ref{vol}).
  \item[ii)] Alternatively, by making the transformation $\widetilde{\mathbf{F}} = \mathbf{P}
    \mathbf{D}_{\widetilde{\mathbf{F}}} \mathbf{P}'$ in (\ref{dfbetaII2}) and integrating over
    $ \mathbf{P} \in \mathcal{O}(r)$ using (\ref{vol}), the proof is completed.
  \item[iii)] A third proof is derived immediately from Theorem \ref{teo1} observing that $\kappa_{i} =
    \delta_{i}^{1/2}$, $i = 1, \dots, r$ and so
    $$\hspace{3.5cm}
      \left(\bigwedge_{i=1}^{r} d\kappa_{i}\right) = 2^{-r}\prod_{i=1}^{r} \delta_{i}^{-1/2}
        \left(\bigwedge_{i=1}^{r} d\delta_{i}\right) \mbox{\hspace{3.5cm}\findemo}
    $$
\end{description}

\begin{rem}\label{rem1}
Finally, from \citet[Lemma 3.1]{m:70} observe that :
\begin{description}
  \item[a)] The density function (\ref{dfbetaII11}) is invariant under any arbitrary matrix
        $\mathbf{\Theta}\in {\mathcal S}_{m}$ and $\mathbf{Y}\in {\mathcal
        L}_{r,m}^{+}(r_{_{\mathbf{\Xi}}}) \sim \mathcal{N}_{m \times r}^{m, r_{_{\mathbf{\Xi}}}}(\mathbf{0},
        \mathbf{\Theta} \otimes \mathbf{\Xi})$.
  \item[b)] Similarly, Theorem \ref{teoBII2} is invariant if $\mathbf{A} \in {\mathcal S}_{m}^{+}(n) \sim
        \mathcal{PW}_{m}^{n}\left(n,\mathbf{\Theta} \right)$, independent of  $\mathbf{Y}\in
        {\mathcal L}_{r,m}^{+}(r) \sim \mathcal{N}_{m \times r}(\mathbf{0}, \mathbf{\Theta}
        \otimes \mathbf{I}_{r})$, where $\mathbf{\Theta}\in {\mathcal S}_{m}$.
  \item[c)] Therefore, Theorem \ref{teochibetaII}, too, is invariant under conditions
        established in items \textbf{a)} and \textbf{b)}.
\end{description}
\end{rem}

\section{Doubly singular beta type I distribution}\label{sec4}

\begin{thm}\label{teoBI1}
Let $\mathbf{A}^{1/2}\mathbf{A}^{1/2} = \mathbf{A} \in {\mathcal S}_{m}^{+}(n) \sim
\mathcal{PW}_{m}^{n}\left(n,\mathbf{I}_{m} \right)$ be independent of  $\mathbf{Y}\in
{\mathcal L}_{r,m}^{+}(r_{_{\mathbf{\Xi}}}) \sim \mathcal{N}_{m \times r}^{m,
r_{_{\mathbf{\Xi}}}}(\mathbf{0}, \mathbf{I}_{m} \otimes \mathbf{\Xi})$, $\mathbf{\Xi}
\in {\mathcal S}_{r}(r_{_{\mathbf{\Xi}}})$ with $ r_{_{\mathbf{\Xi}}} \leq r$ and $m
\geq q \geq r_{_{\mathbf{\Xi}}} > 0$. The doubly singular matrix variate beta type I
distribution is defined as
\begin{eqnarray*}
% \nonumber to remove numbering (before each equation)
   \mathbf{U} &=& \left((\mathbf{A}+\mathbf{Y}\mathbf{\Xi}^{-}\mathbf{Y})^{1/2}\right)^{+}
  \left(\mathbf{Y}\mathbf{\Xi}^{-}\mathbf{Y}\right)
  \left(\mathbf{A}+\mathbf{Y}\mathbf{\Xi}^{-}\mathbf{Y})^{1/2}\right)^{+} \\
     &=& \left(\mathbf{A}+\mathbf{B})^{1/2}\right)^{+} \mathbf{B}
    \left(\mathbf{A}+\mathbf{B})^{1/2}\right)^{+}
\end{eqnarray*}
where $\mathbf{B} =(\mathbf{Y}\mathbf{\Xi}^{-}\mathbf{Y}) \in {\mathcal
S}_{m}^{+}(r_{_{\mathbf{\Xi}}}) \sim \mathcal{PW}_{m}^{r_{_{\mathbf{\Xi}}}}
\left(r_{_{\mathbf{\Xi}}}, \mathbf{I}_{m} \right)$. The density function of
$\mathbf{U}$ is
\begin{equation}\label{dfbetaI1}
     \pi^{-nr_{_{\mathbf{\Xi}}}}
    \frac{ \Gamma_{n+r_{_{\mathbf{\Xi}}}}[(n + r_{_{\mathbf{\Xi}}})/2]}
    {\Gamma_{n}[n/2]\Gamma_{r_{_{\mathbf{\Xi}}}}[r_{_{\mathbf{\Xi}}}/2]}
    |\mathbf{D}_{\mathbf{U}}|^{(r_{_{\mathbf{\Xi}}}-m-1)/2} |\mathbf{I}_{m} -
    \mathbf{U}|^{\left(n-m-1\right)/2} (d\mathbf{U}),
\end{equation}
where $\mathbf{U} = \mathbf{Q}_{1} \mathbf{D}_{\mathbf{U}} \mathbf{Q}'_{1}$ is the
nonsingular part of the spectral decomposition, with $\mathbf{Q}_{1} \in
\mathcal{V}_{r_{_{\mathbf{\Xi}}},m}$ and $\mathbf{D}_{\mathbf{U}} =
\diag(\lambda_{1}, \dots, \lambda_{r_{_{\mathbf{\Xi}}}})$, $1 > \lambda_{1}> \cdots >
\lambda_{r_{_{\mathbf{\Xi}}}}
>0$ and where $(d\mathbf{U})$ denotes the Hausdorff measure.
\end{thm}
\textit{Proof}. From \citet{dggj:06} it is known that if $\mathbf{F} \in
\mathcal{S}^{+}_{m}(r_{_{\mathbf{\Xi}}})$ is a singular matrix variate beta type II
distribution, then the matrix variate $\mathbf{U} = \mathbf{I}_{m} - (\mathbf{I}_{m}
+ \mathbf{F})^{-1}$ has a singular beta type I distribution; moreover
$$
  (d\mathbf{F})= |\mathbf{I}_{r_{_{\mathbf{\Xi}}}} - \mathbf{D}_{\mathbf{U}}|^{-(m+1-r_{_{\mathbf{\Xi}}})/2}
        |\mathbf{I}_{m} - \mathbf{U}|^{-(m+1+r_{_{\mathbf{\Xi}}})/2}(d\mathbf{U})
$$
where $\mathbf{U} = \mathbf{Q}_{1}\mathbf{D}_{\mathbf{U}} \mathbf{Q}'_{1}$ with
$\mathbf{Q}_{1} \in \mathcal{V}_{r_{_{\mathbf{\Xi}}},m}$ and $\mathbf{D}_{\mathbf{U}}
= \diag(\lambda_{1}, \dots \lambda_{r_{_{\mathbf{\Xi}}}})$, $1 > \lambda_{1}> \cdots
> \lambda_{r_{_{\mathbf{\Xi}}}}>0$. The proof follows from making the transformation $\mathbf{F}
=(\mathbf{I}_{m} - \mathbf{U})^{-1}- \mathbf{I}_{m}$ in (\ref{dfbetaII11}), noting
that $\mathbf{F} = (\mathbf{I}_{m} - \mathbf{U})^{-1}\mathbf{U}$ and
$\mathbf{D}_{\mathbf{F}} = (\mathbf{I}_{r_{_{\mathbf{\Xi}}}} -
\mathbf{D}_{\mathbf{U}})^{-1}\mathbf{D}_{\mathbf{U}}$; see \citet{dggj:06}. \findemo

\begin{thm}\label{teoBI2}
Let $\mathbf{A}^{1/2}\mathbf{A}^{1/2} = \mathbf{A} \in {\mathcal S}_{m}^{+}(n) \sim
\mathcal{PW}_{m}^{n}\left(n,\mathbf{I}_{m} \right)$, independent of  $\mathbf{Y}\in
{\mathcal L}_{r,m}^{+}(r) \sim \mathcal{N}_{m \times r}(\mathbf{0}, \mathbf{I}_{m}
\otimes \mathbf{I}_{r})$, and $m \geq n \geq r > 0$. The doubly singular matrix
variate beta type II distribution is defined as
$$
  \widetilde{\mathbf{U}} = \mathbf{Y}'(\mathbf{A} + \mathbf{YY}')^{+} \mathbf{Y},
$$
and its density function is given by
\begin{equation}\label{dfbetaI2}
    \pi^{-r(r+2n-m)/2}
    \frac{ \Gamma_{n+r}[(n + r)/2]} {\Gamma_{n}[n/2]\Gamma_{r}[m/2]}
    |\widetilde{\mathbf{U}}|^{(m- r-1)/2} |\mathbf{I}_{r} -
    \widetilde{\mathbf{U}}|^{-\left(n-m-1\right)/2} (d\widetilde{\mathbf{U}}),
\end{equation}
where $(d\widetilde{\mathbf{U}})$ denotes the Lebesgue measure.
\end{thm}
\textit{Proof}. The proof is obtained immediately by making the transformation
$\widetilde{\mathbf{F}} = (\mathbf{I}_{m} - \widetilde{\mathbf{U}})^{-1}
\widetilde{\mathbf{U}}$ with $(d\widetilde{\mathbf{F}}) = |\mathbf{I}_{m} -
\widetilde{\mathbf{U}}|^{-(r+1)}$ in (\ref{dfbetaII2}). \findemo

\begin{thm}\label{teochibetaI}
Let $\mathbf{U} \in \mathcal{S}_{m}^{+}(r)$ be a matrix variate with density function
(\ref{dfbetaI1}) and $r_{_{\mathbf{\Xi}}} = r$ (or let $\widetilde{\mathbf{U}} \in
\mathcal{S}_{r}$, matrix variate with density function (\ref{dfbetaI2})). Then the
joint density function of $\lambda_{1}, \dots, \lambda_{r}$, the nonzero eigenvalues
of $\mathbf{U}$ (or $\widetilde{\mathbf{U}}$), is
\begin{eqnarray}\label{chibetaI}
    \pi^{-r(2n-m)/2}\frac{\Gamma_{n+r}[(n+r)/2]}{\Gamma_{n}[n/2] \Gamma_{r}[r/2]
    \Gamma_{r}[m/2]}\nonumber\hspace{6cm} \\
    \prod_{i = 1}^{r}\lambda_{i}^{(m-r-1)/2} (1-\lambda_{i})^{(n-m-1)/2}
    \prod_{i<j}\left(\lambda_{i} - \lambda_{j}\right)\left(\bigwedge_{i=1}^{r} d\lambda_{i}\right)
\end{eqnarray}
where $1 > \lambda_{1} > \cdots > \lambda_{r} > 0$.
\end{thm}
\textit{Proof}. The proof can be obtained in any of the following ways.
\begin{description}
  \item[i)] By making the transformation $\mathbf{U} = \mathbf{H}_{1}
    \mathbf{D}_{\mathbf{U}} \mathbf{H}_{1}'$ in (\ref{dfbetaI1}) with $r_{_{\mathbf{\Xi}}} =
    r$ and integrating over $ \mathbf{H}_{1} \in \mathcal{V}_{r,m}$ using
    (\ref{vol}).
  \item[ii)] Alternatively, by making the transformation $\widetilde{\mathbf{U}} = \mathbf{G}
    \mathbf{D}_{\widetilde{\mathbf{U}}} \mathbf{G}'$ in (\ref{dfbetaI2}) and integrating over
    $ \mathbf{G} \in \mathcal{O}(r)$ using (\ref{vol}).
  \item[iii)] A third proof is derived immediately from Theorem \ref{teochibetaII}, observing that
    $\lambda_{i} = (1-\delta_{i})^{-1}\delta_{i}$, $i = 1, \dots, r$, and so
    $$\hspace{3.5cm}
      \left(\bigwedge_{i=1}^{r} d\delta_{i}\right) = \prod_{i=1}^{r} (1-\lambda_{i})^{-2}
        \left(\bigwedge_{i=1}^{r} d\lambda_{i}\right) \mbox{\hspace{3.5cm}\findemo}
    $$
\end{description}

Finally, observe that conclusions analogous to those set out in note 1 can be
obtained for the case of the matrix variate beta type I distribution.

\section{Matricvariate inverted $t$ distribution}\label{sec5}

Let us now find the matricvariate inverted $t$ distribution, also termed the
matricvariate Pearson type II distribution. In this case, we determine its density
function from the doubly singular beta type I distribution (\ref{dfbetaI1}).

\begin{thm}\label{TeoIT}
Let $\mathbf{R} \in {\mathcal L}_{r,m}^{+}(r)$ be the random matrix,
$$
   \mathbf{R} = \left [\left (\mathbf{A} + \mathbf{YY}'\right)^{+}\right]^{1/2} \mathbf{Y},
$$
where $\mathbf{A}^{1/2}\mathbf{A}^{1/2} = \mathbf{A} \in {\mathcal S}_{m}^{+}(n) \sim
\mathcal{PW}_{m}\left(n,\mathbf{I}_{m}\right)$, independent of  $\mathbf{Y}\in
{\mathcal L}_{r,m}^{+}(r) \sim \mathcal{N}_{m \times r}(\mathbf{0}, \mathbf{I}_{m}
\otimes \mathbf{I}_{r})$ and $0 < r \leq n \leq m$. Then the matricvariate
$\mathbf{R}$ has the density function
\begin{equation}\label{mcR}
   dF_{\mathbf{R}}(\mathbf{R}) = \pi^{-r(r+2n)/2}
    \frac{\Gamma_{n+r}[(n+r)/2]}{\Gamma_{n}[n/2]}
    |\mathbf{I}_{m}-\mathbf{RR}'|^{-(n+r)/2}(d\mathbf{R}),
\end{equation}
where $\mathbf{I}_{m}-\mathbf{RR}'> \mathbf{0}$ and $(d\mathbf{R})$ denotes the
Lebesgue measure.
\end{thm}

\textit{Proof}. Observe that $\mathbf{U} = \mathbf{RR}' = \left [\left (\mathbf{A}+
\mathbf{YY}'\right)^{+}\right)]^{1/2}\mathbf{Y}\mathbf{Y}'\left [\left (\mathbf{A}+
\mathbf{YY}'\right)^{+}\right)]^{1/2}$. Then from \citet{dgm:97},
$$
  (d\mathbf{R})= 2^{-r}\left|\mathbf{D}_{\mathbf{U}}\right|^{(r-m-1)/2}(d\mathbf{U})(\mathbf{G}'d\mathbf{G}),
$$
where $\mathbf{R} = \mathbf{H}_{1}\mathbf{D}_{\mathbf{R}}\mathbf{G}'$ and $\mathbf{U}
= \mathbf{RR}' = \mathbf{H}_{1}\mathbf{D}_{\mathbf{R}}\mathbf{G}'\mathbf{G}
\mathbf{D}_{\mathbf{R}} \mathbf{H}'_{1} = \mathbf{H}_{1}\mathbf{D}_{\mathbf{R}}^{2}
\mathbf{H}'_{1} = \mathbf{H}_{1}\mathbf{D}_{\mathbf{U}} \mathbf{H}'_{1}$ with
$\mathbf{D}_{\mathbf{U}} = \mathbf{D}_{\mathbf{R}}^{2}$, $\mathbf{H}_{1} \in
\mathcal{V}_{r,m}$, $\mathbf{G} \in \mathcal{O}(r)$ and $\mathbf{D}_{\mathbf{R}} =
\diag(\tau_{1}, \dots \tau_{r})$, $1 > \tau_{1} > \cdots > \tau_{r}>0$. Thus
$$
  (d\mathbf{U}) =
  2^{r}|\mathbf{N}|^{-(r-m-1)/2}(d\mathbf{R})(\mathbf{G}'d\mathbf{G})^{-1}.
$$
On substituting in (\ref{dfbetaI1}) and integrating over $\mathbf{G} \in O(r)$ using
(\ref{vol}), the proof is complete. \findemo

Alternatively, Theorem \ref{TeoIT} can be proved from the distribution of
$\widetilde{\mathbf{U}}$ in an analogous way.

\section*{Conclusions}

\citet{dggj:08c} studied the density of the singular matricvariate $t$ distribution
and its application to sensitivity analysis. In this work we study some properties of
this distribution, namely the distribution of a linear transformation and the joint
density function of its singular values. \citet{u:94}, \citet{dg:97, dggj:08a}
studied the singular matrix variate beta type I and II distributions in central and
noncentral cases. Now we obtain the doubly singular matrix variate beta type I and II
distributions in the central case and find the joint density function of their
eigenvalues. Finally, we obtain the central density function of the Pearson type II
matricvariate.

\section*{Acknowledgments}

This research work was partially supported  by CONACYT-M\'exico, research grant no.
81512 and IDI-Spain, grants FQM2006-2271 and MTM2008-05785. This paper was written
during J. A. D\'{\i}az- Garc\'{\i}a's stay as a visiting professor at the Department
of Statistics and O. R. of the University of Granada, Spain.

\end{document}